\documentclass{gtart_h}

\def\ifplaintex{\expandafter\ifx\csname documentclass\endcsname\relax}

\def\gtp{{\mathsurround=0pt\it $\cal G\mskip-2mu$eometry \&\ 
$\cal T\!\!$opology $\cal P\!$ublications}}  

\def\Addressesr{\bigskip
{\small \parskip 0pt \leftskip 0pt \rightskip 0pt plus 1fil \def\\{\par}
\sl\theaddress\par
\medskip
\rm Email:\stdspace\tt\theemail\hfill\rm Received:\qua\receiveddate \par}}

\def\recd{{\small Received:\qua\receiveddate\ifx\reviseddate\relax
\else\qquad Revised:\qua\reviseddate\fi\par}} 

\def\lognumber#1{\def\thelognumber{#1}}
\def\volumenumber#1{\def\thevolumenumber{#1}}
\def\volumeyear#1{\def\thevolumeyear{#1}}
\def\papernumber#1{\def\thepapernumber{#1}}
\def\pagenumbers#1#2{\def\startpage{#1}\def\finishpage{#2}}
\def\published#1{\def\publishdate{#1}}

\def\received#1{\def\receiveddate{#1}}

\def\accepted#1{\def\accepteddate{#1}}

\def\asciiaddress#1{\def\theasciiaddress{#1}}

\let\\\par\let\thelognumber\relax\let\thevolumenumber\relax
\let\thepapernumber\relax\let\thevolumeyear\relax\let\startpage\relax
\let\finishpage\relax\let\publishdate\relax\let\receiveddate\relax
\let\reviseddate\relax\let\accepteddate\relax\let\theasciititle\relax
\let\theasciiauthors\relax\let\theasciiaddress\relax
\let\theasciiabstract\relax

\let\theasciiemail\relax

\ifplaintex
\font\logobig=cmssbx10 scaled 3836
\font\logomed=cmssbx10 scaled 2557
\else
\font\logobig=cmssbx10 scaled 4200
\font\logomed=cmssbx10 scaled 2800
\fi

\long\def\makeagttitle{   
\count0=\startpage
\agt\hfill      
\hbox to 45truept{\vbox to 0pt{\vglue -13truept{\logomed A\kern -.37em{\logobig 
T}\kern -.38em G}\vss}\hss}
\break
{\small Volume \thevolumenumber\ (\thevolumeyear)
\startpage--\finishpage\nl
Published: \publishdate}

\vglue .25truein

{\parskip=0pt\leftskip 0pt plus
1fil\def\\{\par\smallskip}{\Large\bf\thetitle}\par\medskip} \vglue
0.05truein

{\parskip=0pt\leftskip 0pt plus 1fil\def\\{\par}{\sc\theauthors}
\par\medskip}%
 
\vglue 0.03truein 

{\small\leftskip 25truept\rightskip 25truept{\bf Abstract}\stdspace\theabstract

{\bf AMS Classification}\stdspace\theprimaryclass
\ifx\thesecondaryclass\relax\else; \thesecondaryclass\fi\par
{\bf Keywords}\stdspace \thekeywords\par}\vglue 7truept

}   

\ifplaintex
\font\phead=cmsl9 scaled 950
\font\pnum=cmbx10 scaled 913
\font\pfoot=cmsl9 scaled 950
\headline{\vbox to 0pt{\vskip -4.5mm\line{\small\phead\ifnum
\count0=\startpage ISSN 1472-2739 (on-line) 1472-2747 (printed)
\hfill {\pnum\folio}\else\ifodd\count0\def\\{ }%
\ifx\theshorttitle\relax\thetitle\else\theshorttitle\fi\hfill{\pnum\folio}
\else\def\\{ and }{\pnum\folio}\hfill\ifx\theshortauthors\relax\theauthors
\else\theshortauthors\fi\fi\fi}\vss}}
\footline{\vbox to 0pt{\vglue 0mm\line{\small\pfoot\ifnum\count0=\startpage
\copyright\ \gtp\hfill\else
\agt, Volume \thevolumenumber\ (\thevolumeyear)\hfill\fi}\vss}}
\else
\font=cmsl9 scaled 1050
\font\lnum=cmbx10 
\font=cmsl9 scaled 1050
\makeatletter
\def\@oddhead{{\small\lhead\ifnum\count0=\startpage ISSN 1472-2739 
(on-line) 1472-2747 (printed)\hfill {\lnum\number\count0}\else\ifodd\count0
\def\\{ }\ifx\theshorttitle\relax \thetitle \else\theshorttitle\fi\hfill
{\lnum\number\count0}\else\def\\{ and }{\lnum\number\count0}
\hfill\ifx\theshortauthors\relax 
\theauthors\else\theshortauthors\fi\fi\fi}}\def\@evenhead{\@oddhead}
\def\@oddfoot{\small\lfoot\ifnum\count0=\startpage\copyright\ \gtp\hfill\else
\agt, Volume \thevolumenumber\ (\thevolumeyear)\hfill\fi}
\def\@evenfoot{\@oddfoot}
\makeatother
\fi
\let\maketitlepage\makeagttitle
\let\makeshorttitle\maketitlepage
\let\maketitle\maketitlepage

\newwrite\gtoutfile
\long\gdef\makeheadfile{  
{\def\\{, }\def\s{ }
\immediate\openout\gtoutfile head.xxx
\immediate\write\gtoutfile{Proxy-for: \ifx\theasciiauthors\relax
\theauthors\else\theasciiauthors\fi\s<\ifx\theasciiemail\relax\theemail\else\theasciiemail\fi>}
\immediate\write\gtoutfile{\noexpand\\}
\immediate\write\gtoutfile{Authors: \ifx\theasciiauthors\relax
\theauthors\else\theasciiauthors\fi}
{\def\\{ }\immediate\write\gtoutfile{Title: \ifx\theasciititle\relax
\thetitle\else\theasciititle\fi}}
\immediate\write\gtoutfile{Subj-class: GT or SG, GR etc}
\immediate\write\gtoutfile{MSC-class: \theprimaryclass\ifx\thesecondaryclass\relax\else, \thesecondaryclass\fi}
\immediate\write\gtoutfile{Journal-ref: Algebraic and Geometric Topology \thevolumenumber\s
(\thevolumeyear) \startpage-\finishpage}
\immediate\write\gtoutfile{Comments: Published by Algebraic and
Geometric Topology at}
\immediate\write\gtoutfile{\s\s\s  http://www.maths.warwick.ac.uk/agt/AGTVol\thevolumenumber/agt-\thevolumenumber-\thepapernumber.abs.html}
\immediate\write\gtoutfile{\noexpand\\}
\immediate\write\gtoutfile{}
\ifx\theasciiabstract\relax
\immediate\write\gtoutfile{\theabstract}\else
\immediate\write\gtoutfile{\theasciiabstract}\fi
\immediate\write\gtoutfile{}
\immediate\write\gtoutfile{\noexpand\\}
\immediate\write\gtoutfile{}
\immediate\closeout\gtoutfile}}  

\def\maketitlepage{\makeagttitle\makeheadfile}
\let\makeshorttitle\maketitlepage
\let\maketitle\maketitlepage

\lognumber{14}
\volumenumber{4}
\volumeyear{2004}
\papernumber{14}
\published{16 April 2004}
\pagenumbers{243}{272}
\received{2 March 2004}
\accepted{15 April 2004}

\usepackage{amssymb} 
\usepackage{rlepsf}
\usepackage{amsmath}

\theoremstyle{definition}
\newtheorem{Def}{Definition}[section]

\newtheorem{Rem}[Def]{Remark}

\theoremstyle{plain}
\newtheorem{Lem}[Def]{Lemma}
\newtheorem{Thm}[Def]{Theorem}

\newtheorem{Prop}[Def]{Proposition}%

\def\R{\mathbb{R}}
\def\C{\mathbb{C}}

\begin{document}

\title{TQFT's and gerbes}
\authors{Roger Picken}
\asciiaddress{Departamento de Matematica and CEMAT - Centro de Matematica e Aplicacoes\\Instituto Superior Tecnico, Av Rovisco Pais\\
1049-001 Lisboa, Portugal}             
\address{Departamento de Matem\'{a}tica and CEMAT - Centro de Matem\'{a}tica e Aplica\c{c}\~{o}es\\Instituto Superior T\'{e}cnico, Av Rovisco Pais\\
1049-001 Lisboa, Portugal}             
\email{rpicken@math.ist.utl.pt}                 

\begin{abstract}  
We generalize the notion of parallel transport along paths for abelian
bundles to parallel transport along surfaces for abelian gerbes using
an embedded Topological Quantum Field Theory (TQFT) approach. We show
both for bundles and gerbes with connection that there is a one-to-one
correspondence between their local description in terms of
locally-defined functions and forms and their non-local description in
terms of a suitable class of embedded TQFT's.
\end{abstract}


\maketitle

\section{Introduction}
\label{Intro}

Gerbes can be regarded as a higher-order version of the geometry of bundles. They originally appeared in the context of algebraic geometry in a paper by Giraud \cite{Gir}, and the concept was further developed in Brylinski's book \cite{Bry}.
Interest in gerbes has been revived recently following a concrete
approach due to Hitchin and Chatterjee~\cite{Hi99}. Gerbes can be understood
both in terms of local geometry, local functions and forms, and in terms of
non-local geometry, holonomies and parallel transports, and these two
viewpoints are equivalent, in a sense made precise by Mackaay and the
author in \cite{MP}, following on from work by Barrett \cite{Bar91} and
Caetano and the author \cite{CP94}. In \cite{MP} holonomies around
spheres and parallel transports along cylinders were considered. The aim of
the present article is to describe the generalization to parallel transports along
general surfaces with boundary.

The conceptual framework for this
generalization is Topological Quantum Field Theory (TQFT), a notion
introduced by Witten~\cite{Wit88}, which was subsequently axiomatized by
Atiyah~\cite{Ati89}, in a very similar manner to Segal's axiomatic
approach to Conformal Field Theory~\cite{Seg89}. The manifolds involved in
the TQFT now come equipped with maps to a target manifold on which the
bundle or gerbe lives. This brings us into the realm of Segal's category
$\cal C$ \cite{SegCFT} and string connection \cite{Seg01}, as well as
Turaev's Homotopy Quantum Field Theory (HQFT) \cite{Tur99} and a related construction by Brightwell and Turner~\cite{BriTur00}, as well as subsequent
developments by Rodrigues~\cite{Rod01},
Bunke, Turner and Willerton~\cite{BunTurWil} and Turner \cite{Turn04}. 
Here we will base our
approach on a general framework for TQFT and related constructions by
Semi\~{a}o and the author~\cite{PS}, which defines a TQFT to be a certain
type of monoidal functor, without using the cobordism approach.

There are a number of interesting connections between the present work and an early approach by Gaw\c{e}dzki \cite{gaw}, further developed by Gaw\c{e}dzki and Reis in \cite{gaw:rei}. In particular equivalent formulae to equations (\ref{ZpT}) and (\ref{ZXT}) for the bundle and gerbe parallel transport appear.

The article is organized as follows. In section 2 we describe some general
preliminaries. In section 3 we define rank-$1$, embedded, $1$-dimensional
TQFT's, show how they can arise from bundles with connection and prove the
equivalence of bundles with connection and a class of these TQFT's. In
section 4 we outline the corresponding procedure for gerbes and rank-$1$,
embedded, $2$-dimensional TQFT's. Section 5 contains some comments.

\rk{Acknowledgements}
Preliminary ideas on the relation between TQFT's and gerbes were presented
at the XXth Workshop on Geometric Methods in Physics, Bialowie\.{z}a, July
1-7, 2001, in the special session on Mathematical Physics at the National
Meeting of the Portuguese Mathematical Society, Coimbra,  February 5-8,
2002, and at the Klausurtagung of the Sonderforschungsbereich 288 -
Differential Geometry and Quantum Physics, St.
Marienthal,  February 25 -  March 1, 2002. I am grateful to the organizers
of these meetings, A. Odzijewicz, A. Strasburger, J. Mour\~{a}o, D. Ferus,
T. Friedrich and R. Schrader for giving me the opportunity to present that
material.

I am grateful to one of the referees for pointing out the connection with Segal's article \cite{Seg68} discussed in Remarks \ref{X_Urem}  and \ref{X_Urem2}.

This work was supported by the program {\em Programa Operacional
``Ci\^{e}ncia, Tecnologia, Inova\c{c}\~{a}o''} (POCTI) of the
{\em Funda\c{c}\~{a}o para a Ci\^{e}ncia e a Tecnologia} (FCT),
cofinanced by the European Community fund FEDER.

\section{Preliminaries}
\label{Prel}

Let us consider smooth, oriented manifolds $B$ of dimension $d\leq 3$,
with or without boundary, partitioned into $d$-dimensional regions by a
finite, embedded $(d+1)$-valent $(d-1)$-graph. Here a $0$-graph is a
collection of points, a $1$-graph is a graph, and a $2$-graph has smoothly
embedded surfaces incident at each vertex. The valency is the number of
$d$-dimensional regions in the neighbourhood of each vertex. Such
partitions can be obtained from triangulations of the manifold by
dualizing. When the manifold has a boundary, the partition of the manifold
should restrict to a partition of the boundary.

Let $M$ be a fixed smooth real manifold of finite dimension, called
the target manifold, with an open cover ${U}=\left\{U_i\ \vert\ i\in
J\right\}$ such that every non-empty $p$-fold intersection
$U_{i_1\ldots i_p}=U_{i_1}\cap\cdots\cap U_{i_p}$, for any $p$, is
contractible. The objects we will be working with are of the form
\[
(Y,T),
\]
where $Y$ is a smooth map from $B$ to $M$, satisfying a ``constant at the
boundary'' condition, to be described later, and $T$ is a labelled
partition of $B$, such that each region of the partition is labelled by an
open set of $M$ containing its image under $Y$. Morphisms between two such
objects $(Y,T)$ and $(Y',T')$ are smooth, orientation-preserving maps
$f:B\rightarrow B'$, such that $Y=Y'\circ f$, 
with no requirements on $T$ or $T'$. 
There is a natural monoidal
structure in this category, which, on objects, is given by:
\[
(Y,T) \sqcup (Y',T')   = (Y\sqcup Y', T \sqcup T'),
\]
where $Y\sqcup Y'$ is the disjoint union of $Y$ and $Y'$ from $B\sqcup B'$
to $M$, and $T \sqcup T'$ is the
obvious labelled partition of $B\sqcup B'$ obtained from $T$ and $T'$.

Our approach in this article is based on a general framework for TQFT and related constructions, which is the subject of \cite{PS}. Since we are working in a very specific context here, and a number of features simplify, it is unnecessary to develop the full formalism, so we will just sketch the general approach to provide some background. The topological category has as its objects essentially pairs consisting of e.g.\ an oriented manifold and its boundary. The boundaries, called subobjects, together with isomorphisms between them, form a category themselves. The morphisms of the topological category are of two types only, namely isomorphisms and gluing morphisms, where the latter are, roughly speaking, morphisms from an object before gluing along one or more boundary components to the object which results after gluing. The topological category is endowed with a monoidal structure (disjoint union) and an endofunctor (change of orientation). Note that the topological category is not set up in the cobordism framework - we refer to the introduction of \cite{PS} for the reasons for preferring an approach closer to Atiyah's original paper \cite{Ati89}.  A TQFT is a functor from the topological category to an algebraic category whose objects are e.g.\ pairs consisting of a finite-dimensional vector space over ${\C}$ and an element of that space. The topological subobject determines the vector space, and the manifold of which it is the boundary determines the element. The algebraic category is endowed with a manifold structure (tensor product) and an endofunctor (passing to the vector space with conjugate scalar multiplication), and the TQFT functor respects these structures in a suitable sense. At the level of morphisms the TQFT sends the isomorphisms between subobjects to certain ``unitary'' algebraic morphisms, namely linear transformations preserving so-called evaluations (resembling inner products). When an isomorphism extends to an isomorphism of topological objects the corresponding linear transformation preserves the elements of the vector spaces involved. Likewise the TQFT functor sends topological gluing morphisms to specific algebraic morphisms involving the evaluations, which again preserve the elements. The requirement that these algebraic morphisms preserve elements implies a set of equations for the elements, which have to be solved for the TQFT functor to exist. In the next sections we will be describing bundle and gerbe parallel transport via operational definitions based on this general TQFT framework.

\section{Embedded TQFT's and bundles}
\label{ETQFTb}
\begin{figure}[ht!]
\centerline{\relabelbox\small 
\epsfysize 3.5cm
\epsfbox{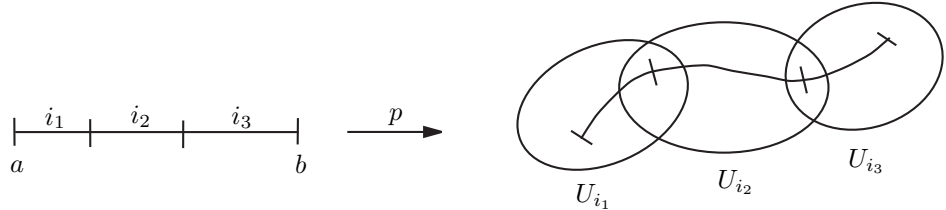} 
\relabel {a}{$a$}
\relabel {b}{$b$}
\relabel {p}{$p$}
\relabel {i1}{$i_1$}
\relabel {i2}{$i_2$}
\relabel {i3}{$i_3$}
\relabel {Ui1}{$U_{i_1}$}
\relabel {Ui2}{$U_{i_2}$}
\relabel {Ui3}{$U_{i_3}$}
\endrelabelbox }
\caption{A path $(p,T)$}
\label{path}
\end{figure}

We will concentrate, for the time being, on $0$ and $1$-dimensional
objects, as described in the previous section. We will call the
$0$-dimensional objects subobjects, and write irreducible subobjects as
\[
(y^\pm , i),
\]
where $y^\pm$ denotes the map from the positively or negatively oriented
point to $M$, whose image is $y\in M$, and $i$ labels an open set of $M$
containing $y$. General subobjects are finite disjoint unions of irreducible
subobjects. An irreducible $1$-dimensional object is of the form
\[
(p,T),
\]
where $p:[a,b]\rightarrow M$ is a smooth path, constant in $[a,
a+\epsilon[$ and $]b-\epsilon, b]$, for some $\epsilon >0$, and $T$ is
given by a set of points $a=x_0<x_1< \cdots < x_N=b$, with a labelling of
$e_\alpha = [x_{\alpha -1}, x_\alpha]$ by $i_\alpha\in J$, such that
$p(e_\alpha)\subset U_{i_\alpha}$ (see Figure \ref{path}). The interval
$[a,b]$ is taken to be oriented in the direction from $a$ to $b$ and the
boundary subobject of $(p,T)$ is $(p(a)^-,i_1)\sqcup (p(b)^+,i_N)$.

As discussed at the end of the previous section there are three types of morphism which enter the description of the topological category: isomorphisms between subobjects, isomorphisms between objects and
gluing morphisms between objects. For the present situation, isomorphisms
between (irreducible) subobjects are described by
\[
(y^\pm, i) \stackrel{\rm id}{\rightarrow} (y^\pm, j)
\]
which will be denoted
\[
(y^\pm, i,j).
\]
Isomorphisms and gluing morphisms between objects are generated (via
composition and disjoint union) by basic isomorphisms
\[
(p,T)\stackrel{f}{\rightarrow}(p',T'),
\]
where $f:[a,b]\stackrel{\cong}{\rightarrow} [a',b']$ is an
orientation-preserving diffeomorphism, and basic gluing morphisms
\[
(p,T)\sqcup (p',T')\stackrel{f}{\rightarrow} (p\circ p',T\circ T')
\]
where $f:[a,b]\sqcup [b,c]\rightarrow [a,c]$ is defined by
$f(x)=x$ on both $[a,b]$ and $[b,c]$,
$p\circ p'$ is given by $(p\circ p')(x)=p(x)$ for
$x\in [a,b]$ and  $(p\circ p')(x)=p'(x)$ for
$x\in [b,c]$, and $T\circ T'$ is the natural labelled partition of $[a,c]$
obtained from $T$ and $T'$. Note that, for this morphism to exist, we 
must have  $p(b)=p'(b)$, since $f(b)=b$ for $b$ in either interval. 
Also note that
$p\circ p'$ is smooth because of the constant condition at the
endpoints of the paths. This method of handling products of smooth loops
and paths was introduced in \cite{CP94}.

\begin{Rem}
The subobjects $\{ (y^\pm, i)\}_{y\in M, \,\, i\in J}$ and the isomorphisms between them form a category, which we will call $M^0_{U}$, where $ U$ denotes the cover $\left\{ U_i \right\}_{i\in J}$ and the superscript indicates the dimension of the subobjects. This category is similar to the category $X_U$ introduced by Segal in \cite{Seg68} for a cover 
$U=\left\{ U_i \right\}_{i\in J}$ of a space $X$ by subsets of $X$. In 
$X_U$ the objects are of the form $(y, U_\sigma)$ where $\sigma$ is a finite subset of $J$, $U_\sigma$ denotes $\cap_{i \in \sigma} U_i$ and $y\in U_\sigma$. The orientation of points is disregarded in $X_U$. The morphisms of $X_U$ may be characterized as follows: for every inclusion $\sigma\subset\tau$ and for every $y\in U_\tau$, there is a unique morphism, which we will denote $(y,\tau,\sigma)$, from $(y, U_\tau)$ to $(y, U_\sigma)$. Shortly we will return to a comparison of $M^0_{U}$ and $X_U$ in the context of bundles. 
\label{X_Urem}
\end{Rem}

The full definition of a TQFT functor in \cite{PS} can be replaced here by a considerably reduced
operational definition, since in the present situation we are able to make a number of simplifying choices. We take all vector spaces in the algebraic category to be of dimension $1$, which allows the ``unitary'' algebraic morphisms corresponding to topological isomorphisms of subobjects to be replaced by complex numbers of norm $1$. Likewise we need only write down the equations which arise from applying the TQFT functor to basic isomophisms and gluing morphisms, omitting equations relating to the monoidal structure (disjoint union, empty set), which are understood. The elements of the $1$-dimensional vector spaces corresponding to topological objects $(p,T)$ may again be identified with complex numbers of norm $1$. Thus we are led to the following definition.
\begin{Def}
A rank-$1$, embedded, $1$-dimensional TQFT is a pair of assignments
$(Z',Z)$ 
\begin{eqnarray*}
(y^\pm, i,j) & \mapsto &  Z'(y^\pm, i,j) \in U(1) \\
(p,T)  & \mapsto & Z(p,T) \in U(1)
\end{eqnarray*}
such that
\begin{itemize}
\item[i)]  for isomorphisms between subobjects
\[
 Z'(y^\pm, i,j)  Z'(y^\pm, j, k) =  Z'(y^\pm, i,k),
\]
\item[ii)]  for basic isomorphisms $(p,T)\stackrel{f}{\rightarrow}(p',T')$
\[
Z(p',T')= Z'(p(a)^-, i_1, i'_1)Z(p,T) Z'(p(b)^+, i_N, i'_{N'}),
\]
\item[iii)] for basic gluing morphisms $(p,T)\sqcup
(p',T')\stackrel{f}{\rightarrow} (p\circ p',T\circ T')$
\[
Z(p\circ p',T\circ T') = Z(p,T) Z'(p(b)^+, i_N, i'_1)   Z(p',T').
\]
\end{itemize}
\label{def1dTQFT}
\end{Def}

The terminology rank-$1$ is the same as that used for Homotopy Quantum Field Theories in \cite{BunTurWil}.
We remark that an analogous definition can be given of rank-$n$ embedded
TQFT's by replacing $U(1)$ by $U(n)$. We will comment on this in section
\ref{Comm}. 
Next we state some simple properties of TQFT's which will be needed later on.
\begin{Prop}
For a rank-$1$, embedded, $1$-dimensional TQFT $(Z',Z)$ the following properties hold:
\begin{itemize}
\item[\rm a)] $Z'(y^\pm, i,i)=1$,
\item[\rm b)] for any interval $[a,b]$, 
\[
Z(t_{y[a,b]},i) =1
\]
where $ t_{y[a,b]}$ denotes the constant map from $[a,b]$ to $y$, and $i$ denotes the labelled partition which assigns $i$ to the whole interval $[a,b]$,
\item[\rm c)] $Z(t_y,ij)=Z'(y^+,i,j)$, where $ t_{y}$ denotes the constant map from the standard interval $[0,1]$ to $y$, and
$ij$ is the labelled partition which labels $[0,1/2]$ with $i$ and
$[1/2,1]$ with $j$,
\item[\rm d)] $Z'(y^-,i,j)= Z'(y^+,i,j)^{-1}$.

\end{itemize}
\label{prop}
\end{Prop}

\proof a) is obvious, from i) of Definition \ref{def1dTQFT}. To show b) we first use ii) applied to an obvious isomorphism to obtain
\[
Z(t_{y[a,b]},i)= Z(t_{y[c,d]},i)
\]
for any intervals $[a,b]$ and  $[c,d]$, and then derive
\[
\begin{array}{lcl}
Z(t_{y[a,b]},i)& \stackrel{ii)}{=} & Z(t_{y[a,c]}\circ t_{y[c,b]},i \circ i)
\\
& \stackrel{iii)\,a)}{=} & Z(t_{y[a,c]},i) Z(t_{y[c,b]},i).
\end{array}
\]
c) follows from applying iii) to the obvious gluing morphism from $(t_{y[0,1/2]},i) \sqcup (t_{y[1/2,1]},j)$ to $(t_{y},ij)$ and using b). To show d), we apply ii) to any trivial path as follows:
\[
Z(t_y,j)=Z'(y^-,i,j) Z(t_y,i) Z'(y^+,i,j)
\]
and use b). \endproof

Suppose we are given a $U(1)$-bundle with connection on the target manifold
$M$, i.e.\ a collection of (smooth) transition functions
\[
g_{ij}: U_{ij}\rightarrow U(1)
\]
and (smooth) local connection $1$-forms
\[
A_j \in \Lambda ^1(U_j),
\]
satisfying
\begin{itemize}
\item[B1)] $g_{ij}g_{jk}=g_{ik}$ on $U_{ijk}$,
\item[B2)] $i(A_k-A_j)=d\log g_{jk}$ on $U_{jk}$.
\end{itemize}
B1) is called the cocycle condition. It implies also the equation
$g_{ji} = g_{ij}^{-1}$ on $U_{ij}$. 
The curvature of this connection is
the globally-defined $2$-form $F\in \Lambda ^2(M)$, given locally by
\begin{equation}
F=dA_i.
\label{FdA}
\end{equation}

\begin{Thm}
A $U(1)$-bundle with connection on $M$ gives rise to a  rank-$1$, embedded,
$1$-dimensional TQFT with target $M$ via
\begin{equation}
Z'(y^+,i,j) = g_{ij}(y), \quad Z'(y^-,i,j) = \bar{g}_{ij}(y)
\label{Z'y}
\end{equation}
\begin{equation}
\begin{array}{lll}
Z(p,T)&=& \exp i\int_{x_0}^{x_1} p^\ast(A_{i_1}). g_{i_1i_2}(p(x_1)).
\exp i\int_{x_1}^{x_2} p^\ast(A_{i_2}).  g_{i_2i_3}(p(x_2)) \dots
\\[15pt]
& & \exp i\int_{x_{N-1}}^{x_N} p^\ast(A_{i_N}).
\label{ZpT}
\end{array}
\end{equation}
\label{thmbundletqft}
\end{Thm}
\proof i) of Definition \ref{def1dTQFT} follows immediately from the
cocycle condition B1). The proof of ii) is deferred until after Lemma 
\ref{lemma2stokes}. 
iii) is immediate.
\medskip

Let us now introduce the following shorthand notation for the
right-hand-side of equation (\ref{ZpT}):
\[
\exp i \int_{(p,T)} (g,A)
\]
and an analogous notation for $(\gamma , T)$, where $\gamma:S^1\rightarrow
M$ is a closed curve, and $T$ is a labelled partition of $S^1$. We will
also now introduce $2$-dimensional objects $(H,T)$, where $H$ is a smooth map
from an oriented surface $S$ with or without boundary, and $T$ is a
labelled partition of $S$,  by means of an embedded $3$-valent graph, into
faces $f_\alpha$, such that $H(f_\alpha) \subset U_{i_\alpha}$. When $S$
has a boundary, we will denote by $(\partial H, \partial T)$ the
$1$-dimensional object obtained by restricting $H$ and the labelled
partition $T$ to the oriented boundary $\partial S$ of $S$.

We now have the following useful Lemma:

\begin{Lem} For a $2$-dimensional object $(H,T)$,
\[
\exp i \int_{(\partial H, \partial T)} (g,A) = \exp i \int_S H^\ast (F).
\]
\label{lemma2stokes}
\end{Lem}
\proof We have
\[
\exp i \int_{(\partial H, \partial T)} (g,A) =
\prod_\alpha \exp i \int_{\partial f_\alpha} H^\ast(A_{i_\alpha}),
\]
since the transition function contributions around the boundary may
be replaced by the factors on the right-hand-side of the above equation
corresponding to the internal
edges of the partition, using the relation for an edge $e$ of the
partitioning graph from $x$ to $y$ in $S$, between faces labelled $j$ and
$k$:
\[
\exp i \int_e H^\ast(A_k-A_j)=\exp  \int_{\partial e} H^\ast(\log
g_{jk}) = g^{-1}_{jk}(H(x)) g_{jk}(H(y)),
\]
and the fact that the
contributions of the transition functions cancel at each internal vertex
because of the cocycle condition. Using equation (\ref{FdA}) and Stokes'
theorem the result follows. \endproof

We now complete the proof of Theorem \ref{thmbundletqft}.

\begin{figure}[ht!]
\centerline{\relabelbox\small 
\epsfysize 8cm
\epsfbox{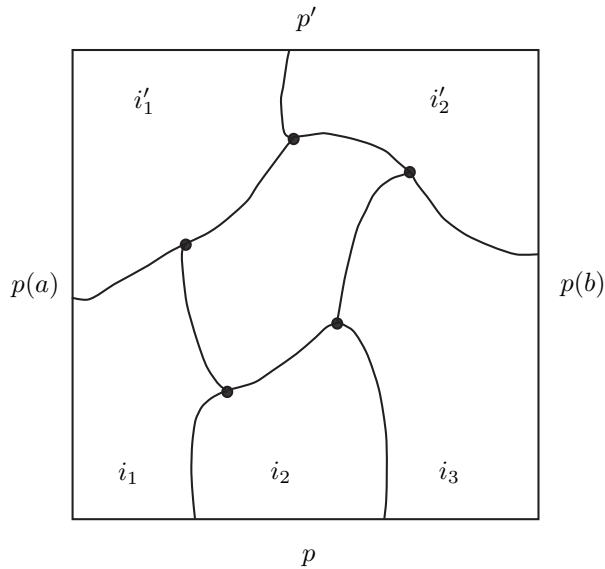} 
\relabel {p}{$p$}
\relabel {p1}{$p'$}
\relabel {pa}{$p(a)$}
\relabel {pb}{$p(b)$}
\relabel {i2}{$i_2$}
\relabel {i1}{$i_1$}
\relabel {i3}{$i_3$}
\relabel {j1}{$i'_1$}
\relabel {j2}{$i'_2$}
\endrelabelbox }
\caption{Proof of Theorem \ref{thmbundletqft}, part ii)}
\label{proofii}
\end{figure}

\noindent\proof  To show  ii) of Definition \ref{def1dTQFT}, we introduce a
homotopy $(H,\tilde{T})$ from $p$ to $p'$, whose image is contained in the
shared image of   $p$ and $p'$, and such that $\partial \tilde{T} $
coincides with $T$ and $T'$ on the top and bottom edges, and changes from
$i_1$ to $i'_1$ and from $i_N$ to $i'_{N'}$ on the two sides (see Figure
\ref{proofii}). It is easy to see that such a partition $\tilde{T}$ exists,
by constructing partitions for the elementary moves of splitting an
interval labelled with $i$ into two regions labelled with $i$ (and the
inverse move), and relabelling a single interval. Applying the Lemma and
using $H^\ast(F)=0$, we get \[ g_{i_1i'_1}(p(a)) Z(p,T)
g_{i_Ni_{N'}'}(p(b)) Z^{-1}(p',T')=1, \] which implies ii). Here we have
used the equation
\[ Z(p,T) Z(p^{-1}, T^{-1})=1, \] where $p^{-1}:[a,b]\rightarrow M$ is
given by $p^{-1}(x)=p(a+b-x)$ and $T^{-1}$ is the corresponding inverse
labelled partition, which follows from applying the Lemma to the identity
homotopy from $(p,T)$ to $(p,T)$. \endproof

\begin{Rem}
The assignment $Z'$ of Theorem \ref{thmbundletqft} defines a functor from the category $M_U^0$, introduced in Remark \ref{X_Urem}, to $U(1)$, regarded as a $1$-object category, because of equation i) of Definition \ref{def1dTQFT}. Conversely it is clear that, subject to a smoothness restriction, any such functor defines transition functions, so that we have a correspondence between smooth functors from $M_U^0$ to $U(1)$ and smooth bundles on $M$ defined by transition functions $g_{ij}$ on $U_{ij}$. In what follows we will be extending this correspondence to one between a suitable class of TQFT's as in Definition \ref{def1dTQFT} and smooth bundles with connection. It is interesting to contrast this with Segal's functorial approach to principal bundles in \cite{Seg68}. A functor from the category $X_U$ (discussed in Remark \ref{X_Urem}) to a topological group $G$ implies choosing transition functions on the cover $U'$ of $X$, whose elements are all finite intersections of elements of the cover $U$, via $(y,\tau,\sigma) \mapsto g_{\sigma\tau}(y) \in G$. Transition functions on the cover $U$ are obtained by
\[
g_{ij}(y) = g_{\{i\}\{i,j\}}(y) (g_{\{j\}\{i,j\}}(y))^{-1}.
\]
Then the statement is that (smooth) principal $G$-bundles over $X$ defined with respect to the cover $U$ correspond to equivalence classes of (smooth) functors from $X_U$ to $G$.
\label{X_Urem2}
\end{Rem}

 The TQFT's obtained in this way from bundles with connection satisfy two
special properties.

\begin{Thm} The rank-$1$, embedded, $1$-dimensional TQFT's of Theorem
\ref{thmbundletqft} satisfy
\begin{itemize}
\item[\rm 1)] thin invariance: if $(p,T)$ and $(p',T')$ are such
that $i_1=i'_1$,
$i_N=i'_{N'}$ and there is a relative homotopy $(H, \tilde{T})$ from $(p,T)$
to $(p',T')$ satisfying ${\rm rank}\, DH \leq 1$ (a so-called thin homotopy
- see \cite{CP94,MP}), then \[ Z(p,T)=Z(p',T'), \] \item[\rm 2)] smoothness:
for any smooth $k$-dimensional family of objects $(p(u), T(u))$, $u\in
U\subset {\R}^k$, $Z(p(u),T(u))$ depends smoothly on $u$. 
\end{itemize}
\label{thinthm} 
\end{Thm} \proof 1) follows directly from the Lemma,
noting that $H^\ast(F)=0$, since ${\rm rank}\, DH\leq 1$. 2) follows from
the fact that the formula for $Z(p,T)$ is smooth, since it involves smooth
functions and integrals. \endproof

We remark that, as pointed out by Bunke et al \cite{BunTurWil}, the thin
invariance also holds if one replaces the homotopy $H$ by a
cobordism which is thin in the same sense as above.

The properties of the previous theorem in fact characterize the TQFT's
obtained from bundles with connection, because of the following theorem.
\begin{Thm}
There is a one-to-one correspondence between bundles with connection on $M$
and smooth, thin invariant, rank-$1$, embedded, $1$-dimensional TQFT's with
target $M$.
\label{btqftcorr}
\end{Thm}
\proof We have already seen the correspondence in one direction (Theorem 
\ref{thmbundletqft} and Theorem \ref{thinthm}).
Suppose we are given an embedded TQFT with the stated properties. We may
reconstruct a bundle with connection on $M$ as follows.
\begin{equation}
g_{ij}(y) = Z(t_y, ij)
\label{gfromZ}
\end{equation}
where $t_y$ is the trivial, constant path from $[0,1]$ to $y\in M$, and
$ij$ is the labelled partition which labels $[0,1/2]$ with $i$ and
$[1/2,1]$ with $j$.
\begin{equation}
i(A_j)_y(v)= \frac{d}{dt} \log Z(q_{t}, j) |_{t=0}
\label{AfromZ}
\end{equation}
where we take a path $q:(-\epsilon , \epsilon)\rightarrow U_i$ such that
$q(0)=y, \, \dot{q}(0)=v$, and set $q_{t}$ to be a path from $y$ to
$q(t)$ along $q$, reparametrized to be constant at the endpoints.
Furthermore $j$ denotes the labelled partition assigning $j$ to the whole
domain of $q_{t}$.

To show that $A_j$ is well-defined, we note first that, given a choice of $q$, the definition does not depend on the choice of $q_t$, because of the fact that the different choices are all thin equivalent, and using the thin invariance of $Z$. Thus it remains to show that the definition does not depend on the choice of $q$. We use the fact that $M$ is a manifold, i.e.\ locally diffeomorphic to an open neighbourhood of the origin in ${\R}^d$, to fix a smooth family of radial paths, denoted $r_x$, from any point $x$ in a sufficiently small neighbourhood of $y$ to $y$ (reparametrized to be constant at the endpoints). 
Now
\begin{equation}
\frac{d}{dt} \log Z(q_{t}, j) |_{t=0} = - \frac{d}{dt} \log Z(r_{q(t)}, j) |_{t=0}
\label{Awelldef}
\end{equation}
since 
\[
\frac{d}{dt} \log Z(q_{t}\circ r_{q(t)}, j) |_{t=0} =0,
\]
where we use Barrett's lemma \cite{Bar91} (see also \cite{CP94, MP}) stating that the trivial loop is a critical point for any holonomy. Note that $Z(.,j)$ restricted to loops based at $y$ contained in $U_j$ satisfies the conditions to be a holonomy in the sense of Barrett's lemma, because of the thin invariance and smoothness of $Z$ (Theorem \ref{thinthm}). Now the right-hand-side of equation (\ref{Awelldef}) is equal to $(dh)_y(v)$, where $h$ is a smooth function defined in a neighbourhood of $y$ by 
$h(x) = \log Z(r_x,j)$, and thus is independent of the choice of $q$ with $\dot{q}(0)=v $.

Now property B1) for the transition functions (\ref{gfromZ}) is immediate, using 
Proposition \ref{prop} c) and i) of Definition \ref{def1dTQFT}. To show B2) consider $\phi : (-\epsilon, \epsilon) \rightarrow i {\R}$ defined by
\[
\phi(t) = \log Z(q_{t} \circ q^{-1}_{t}, j \circ k).
\]
By the thin invariance of $Z$, $\phi$ is constant, since there are natural thin homotopies between the arguments of $Z$ for different values of $t$. Thus
\[
\begin{array}{lcl}
 0 & \stackrel{iii)}{=} & \frac{d}{dt} 
\left[\log Z(q_{t}, j) + \log Z'(q(t),jk)
+ \log Z(q^{-1}_{t}, k) \right] |_{t=0} \\
& = & i(A_j)_y(v) + (d\log g_{jk})_y(v) - i(A_k)_y(v),
\end{array}
\]
which proves B2).

Finally we show the correspondence is one-to-one. Let 
$(g_{ij}, A_j)_{i,j\in J}$ be a bundle with connection, and let $(Z',Z)$ be the TQFT which arises from it (Theorem \ref{thmbundletqft}). Now take 
$(\tilde{g}_{ij},\tilde{A}_j) _{i,j\in J}$ to be the bundle and connection reconstructed from $(Z',Z)$ via equations (\ref{gfromZ}) and (\ref{AfromZ}). We have
\[
\tilde{g}_{ij}(y)=Z(t_y,ij)=Z'(y^+,i,j)=g_{ij}(y)
\]
using Proposition \ref{prop} c) in the second equality, and
\begin{eqnarray*}
i(\tilde{A}_{j})_y(v) &=&\frac{d}{dt} \log Z(q_{t}, j) |_{t=0}=
 i\frac{d}{dt}  \int_0^1 q_t^*(A_j) |_{t=0} \\
&=& i\frac{d}{dt}  \int_0^t q^*(A_j) |_{t=0} = i(A_j)_y(v)
\end{eqnarray*}
where, in the first equality, we choose $q_t: [0,1]\rightarrow M$ defined by $q_t(s)=q(st)$, reparametrized to be constant at the endpoints, in the second equality we use equation (\ref{ZpT}), and, in the third equality, we use the reparametrization invariance of the integral and make a substitution. 

Conversely, let $(Z',Z)$ be a TQFT, as specified, and let 
$(g_{ij}, A_j) _{i,j\in J}$ be the bundle with connection reconstructed from it. Now take $(\tilde{Z}', \tilde{Z})$ to be the TQFT that this bundle with connection gives rise to (Theorem \ref{thmbundletqft}). We have
\[
\tilde{Z}'(y^+,i,j)=g_{ij}(y)=Z(t_y,ij)=Z'(y^+,i,j)
\]
and
\[
\tilde{Z}'(y^-,i,j)= \tilde{Z}'(y^+,i,j)^{-1}= {Z}'(y^+,i,j)^{-1} ={Z}'(y^-,i,j),
\]
using Proposition \ref{prop} c), d) and 
\begin{eqnarray*}
\tilde{Z}(p,T) & = & 
\exp i\int_{x_0}^{x_1} p^\ast(A_{i_1}). g_{i_1i_2}(p(x_1)).
\exp i\int_{x_1}^{x_2} p^\ast(A_{i_2}).  g_{i_2i_3}(p(x_2)) \dots \\
& = & Z(p_1,i_1)Z'(p(x_1)^+, i_1,i_2)Z(p_2,i_2) \dots \\
& = & Z(p,T)
\end{eqnarray*}
where, in the second equality, we have written 
$(p,T)=(p_1,i_1)\circ (p_2,i_2) \dots $, 
and used Proposition \ref{prop} c), together with a calculation below, and in the third equality, we have used iii) of Definition \ref{def1dTQFT}. The calculation required is:
\begin{equation}
\exp i \int_a^b p^*(A_j)= Z(p,j)
\label{eintA=Z}
\end{equation}
for a path $p:[a,b]\rightarrow M$, which is shown as follows. First we have
\begin{eqnarray*}
ip^*(A_j) & = & i(A_j)_{p(x)}(\dot{p}(x))\, dx \\
& = & \frac{d}{dt} \log Z(q_t,j)|_{t=0}\, dx \\
& = & \frac{d}{dx} \log Z(p_x,j) dx
\end{eqnarray*}
where, in the second equality, we have introduced $q:(-\epsilon,\epsilon) \rightarrow M$, defined by $q(s)=p(x+s)$, 
and $q_t:[0,1]\rightarrow M$ given by $q_t(s)=q(ts)$, reparametrized to be constant at the endpoints, and in the third equality, we have introduced $p_x:[0,1] \rightarrow M$, defined by $p_x(s)=p(a+s(x-a))$, reparametrized to be constant at the endpoints, and used the relation
\begin{eqnarray*}
\log Z(q_t,j) & = & \log Z(p_x\circ q_t,j) - \log Z(p_x,j)\\
& = & \log Z(p_{x+t},j) - \log Z(p_x,j)
\end{eqnarray*}
in the limit $t\rightarrow 0$, where the first equality follows from iii)
of Definition \ref{def1dTQFT} and the second equality follows from the thin equivalence of $p_x\circ q_t$ and $p_{x+t}$ and the thin invariance of $Z$. Thus we have
\[
\exp i \int_a^b p^*(A_j)= Z(p_b,j) Z(p_a,j)^{-1} = Z(p,j)
\]
where the final equality follows from the thin equivalence of 
$p_b$ and $p$ and the thin invariance of $Z$, as well as from 
Proposition \ref{prop} b), since $p_a$ is a trivial path. This completes the proof that the maps given by equations (\ref{Z'y}) and (\ref{ZpT}) on the one hand and equations (\ref{gfromZ}) and (\ref{AfromZ}) on the other hand are mutual inverses, i.e.\ that the correspondence is one-to-one.
\endproof

\begin{Rem}
In the approach of Mackaay and the author in \cite{MP}, the role of TQFT's was played by holonomies, which assign a group element to each based loop, up to certain conditions such as multiplicativity. The reconstruction of transition functions $g_{ij}$ was ambiguous there, since it 
depended on fixing paths from the basepoint to each point of $M$, and similarly for the reconstruction of the connection $1$-forms $A_i$. Thus the correspondence between bundles with connection and holonomies (Theorem 3.9 of \cite{MP}) was necessarily only up to equivalence. The advantage of the TQFT approach here, based on paths, is that there is a natural path associated to each point of $M$, namely the trivial path at that point, and there are natural families of paths associated to each vector at a point, namely the paths $q_t$ of the above proof. Thus here the correspondence works without having to descend to equivalence classes. Indeed a TQFT can be regarded as the extension to all paths $(p,T)$ of assignments to certain trivial and infinitesimal paths given by $g_{ij}$ and $A_i$. 
\label{corresp_rem}
\end{Rem}

\section{Embedded TQFT's and gerbes}
\label{ETQFTg}

\begin{figure}[ht!]
\centerline{\epsfysize 4cm\epsfbox{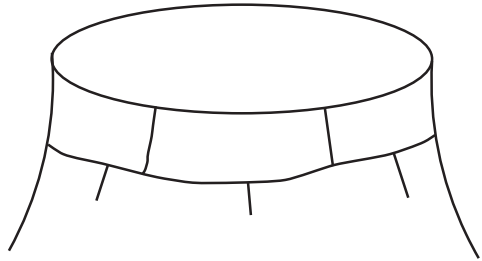} }
\caption{Partition near a boundary circle}
\label{boundpartn}
\end{figure}

We will now sketch the analogous picture for gerbes.
The irreducible subobjects are now $1$-dimensional, and
of the form
\[
(\ell^\pm, T)
\]
where $\ell^\pm : S^1_\pm \rightarrow M$ is a smooth map from the $+$
(anticlockwise) or $-$ (clockwise) oriented circle to $M$, and $T$ is a
labelled partition of $S^1$ given by a choice of points ${a}_1, \dots ,
{a}_N $ on $S^1$, ordered cyclically, together with an assignment to
the edge (arc) $e_\alpha$, between ${a} _{\alpha - 1}$ and ${a}
_\alpha$, of $i_\alpha \in J$, such that $\ell(e_\alpha )\subset
U_{i_\alpha}$.
General subobjects are
finite disjoint unions of irreducible subobjects. An irreducible
$2$-dimensional object is of the form
\[ (X,T), \]
where $X$ maps from an oriented surface with or without boundary $S$ to
$M$, satisfying the following ``constant at the boundary'' condition: at
any boundary circle of $S$ with local coordinates ${a}$ around the
circle and $r$ transversal to the circle, with the boundary circle itself
being $r=r_0$, $X$ is constant in $r$ for a neighbourhood of $r_0$. The
labelled partition $T$ of $S$ is given by a $3$-valent embedded graph in
$S$, and a labelling of the regions (faces) $f_\alpha$ by $i_\alpha\in J$,
such that $X(f_\alpha )\subset U_{i_\alpha}$, and such that the partition
looks like Figure \ref{boundpartn} within the constant neighbourhood of
any boundary circle.

\begin{figure}[ht!]
\centerline{\epsfysize 6cm\epsfbox{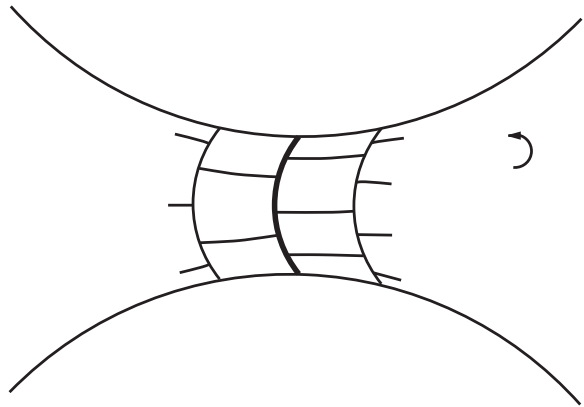} }
\caption{Partition in the neighbourhood of a cut}
\label{cut}
\end{figure}

As before we will consider
isomorphisms between subobjects, basic isomorphisms between objects and
basic gluing morphisms between objects, and will also introduce a fourth type 
of morphism, called partial gluing morphisms. 
(This last type of morphism was not considered in \cite{PS}, since it belongs
 more to the realm of extended TQFT, but we will see that it fits into 
the picture very naturally here.)
Isomorphisms
between (irreducible) subobjects are described by
\[
(\ell^\pm, T) \stackrel{\rm id}{\rightarrow} (\ell^\pm, T')
\]
which will be denoted
\[
(\ell^\pm, T, T').
\]
Basic isomorphisms between $2$-dimensional objects are given by
\[
(X,T)\stackrel{f}{\rightarrow}(X',T'),
\]
where $f:S\stackrel{\cong}{\rightarrow} S'$ is
orientation-preserving, and restricts to the identity at the boundary
circles, i.e.\ $f({a} , r_0)=({a} , r_0')$ with respect to the local
coordinates introduced in the neighbourhood of the boundary circles.
Next, basic gluing morphisms between $2$-dimensional objects are given by
\[
(X,T)\stackrel{f}{\rightarrow} (X', T')
\]
where $X:S\rightarrow M,\, X':S'\rightarrow M$, $S$ is obtained from $S'$
by cutting along a circle $C$ in $S'$, which becomes two boundary circles
in $S$, and $f$ is the identity map everywhere, except at the cut where $f$
maps the pair of boundary circles arising from the cut to $C$ in $S'$. The
partition $T'$ in a neighbourhood of $C$ should be as in Figure
\ref{cut}, with the labelled partitions in the constant regions to the
left and right of $C$ denoted $T'_L,\, T'_R$ respectively. The labelled
partitions $T$ and $T'$ are identical except from the different
interpretations at the cut. Note that we have not written the domain of the
gluing morphism as a disjoint union, since there is the possibility of
self-gluing, a point which is emphasized in \cite{PS}. Finally, a partial gluing morphism is given in the same way as a gluing morphism, except that the cut is along a curve in $S'$ which intersects the boundary of $S'$ transversally, and the labelled partitions on both sides of the curve are identical. 
Here we will assume that $S'$ is divided into two pieces by the cut, i.e.\ 
$S'=S_1\sqcup S_2$, and $(X,T)$ also splits into two objects 
$(X_1:S_1\rightarrow M,T_1)$ and $(X_2:S_2\rightarrow M,T_2)$.
An example of a partial gluing morphism is when 
$S_1$ and $S_2$ are triangular surfaces and are glued along one edge to form a square surface $S'$. Such examples will appear in the proof of Theorem \ref{gtqftcorr}.

\begin{Def}
A rank-$1$, embedded, $2$-dimensional TQFT is a pair of assignments
\begin{eqnarray*}
(\ell ^\pm,T,T') & \mapsto &  Z'(\ell ^\pm,T,T') \in U(1) \\
(X,T)  & \mapsto & Z(X,T) \in U(1)
\end{eqnarray*}
such that
\begin{itemize}
\item[i)]  for isomorphisms between subobjects
\[
 Z'(\ell ^\pm, T, T')  Z'(\ell ^\pm, T', T'') =  Z'(\ell ^\pm, T,T''),
\]
\item[ii)]  for basic isomorphisms $(X,T)\stackrel{f}{\rightarrow}(X',T')$
\[
Z(X',T')= \left( \prod_{ C \in \partial S} Z'(X'|_C^{o(C)},
T|_C,T'|_C) \right) Z(X,T),
\]
where $o(C)=\pm$ is the orientation of $C$ induced from the orientation of $S$,
\item[iii)] for basic gluing morphisms as described above,
corresponding to a cut along $C$ in $S'$
 \[
 Z(X',T') =  Z'(X'|_C, T'_L, T'_R) Z(X,T),
  \]
\item[iv)] for partial gluing morphisms as described above
\[
Z(X',T') = Z(X_1,T_1) Z(X_2,T_2).
\]
\end{itemize}
\label{def2dTQFT}
\end{Def}

\begin{figure}[ht!]
\centerline{\relabelbox\small 
\epsfysize 8cm
\epsfbox{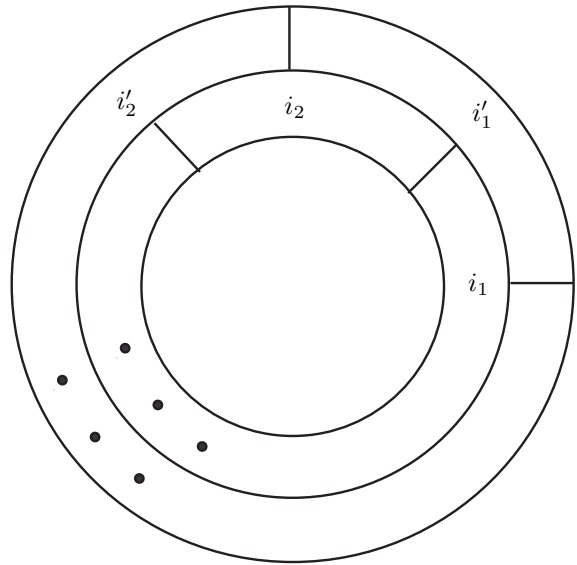} 
\relabel {i1}{$i_1$}
\relabel {i2}{$i_2$}
\relabel {j1}{$i'_1$}
\relabel {j2}{$i'_2$}
\endrelabelbox }
\caption{Partition of the annulus $T\cup T'$}
\label{annulus}
\end{figure}

We remark that the assignment of a phase to a loop with two labelled partitions in i) of the above
definition, suggests a notion of transition functions for bundles on the loop space of the target manifold, and indeed this was made precise in \cite{gaw, gaw:rei}.
We now state some analogous properties to Proposition \ref{prop}.
\begin{Prop}
For a rank-$1$, embedded, $2$-dimensional TQFT $(Z',Z)$ the following properties hold:
\begin{itemize}
\item[\rm a)] $Z'(\ell^\pm, T,T)=1$,
\item[\rm b)] for any $\ell: S^1 \rightarrow M$ and interval $I=[b,c]$ define 
$L_I:S^1\times I \rightarrow M$ by $L_I({a},x)=\ell({a})$. Then
\[
Z(L_I, T) =1,
\]
where $T$ is any partition consisting of regions separated by lines of constant angle $a$,
\item[\rm c)] $Z(L, T\cup T')=Z'(\ell^+, T, T')$, where 
$L:S^1\times [0,1] \rightarrow M$ is given by $L({a},x) = \ell({a})$ and the partition of the annulus $T\cup T'$ is depicted in Figure \ref{annulus} (the inner decomposition is $T$ and the outer is $T'$),
\item[\rm d)] $Z'(\ell^-,T,T')= Z'(\ell^+,T,T')^{-1}$.
\end{itemize}
\label{prop2}
\end{Prop}
\proof a) is immediate, from i) of Definition \ref{def2dTQFT}. To show b), we have $ Z(L_I, T)$ $= Z(L_{I'}, T)$ for any two intervals $I,\,I'$, by Definition
\ref{def2dTQFT} ii) and a). For $I=[b,d]$ and $I'=[d,c]$ we have
\[
Z(L_{I\cup I'}) \stackrel{iii)\,b)}{=} Z(L_I, T) Z(L_{I'}, T)  
\] 
and the result follows. 
c) follows from 
\[
Z(L, T\cup T') \stackrel{iii)}{=} 
Z'(\ell^+, T, T') Z(L_{[0,1/2]}, T) Z(L_{[1/2,1]}, T')\stackrel{b)}{=}
Z'(\ell, T, T').  
\]
To show d), we apply ii) of Definition \ref{def2dTQFT} as follows:
\[
Z(L,T')=Z'(\ell^-,T,T')Z'(\ell^+,T,T') Z(L,T) 
\]
where $L$ is defined in c), and use b). \endproof

Suppose we are given a $U(1)$-gerbe with connection on the target manifold
$M$, i.e.\ a collection of transition functions on triple intersections
\[
g_{ijk}: U_{ijk}\rightarrow U(1)
\]
and local connection $1$- and $2$-forms
\[
A_{jk} \in \Lambda ^1(U_{jk}), \quad F_k \in \Lambda ^2(U_k),
\]
satisfying
\begin{itemize}
\item[G1)] $g_{ijk}=1$ for $i=j$, $i=k$, or $j=k$,
\item[G2)] $g_{ijk}g_{ikl}=g_{jkl}g_{ijl}$ on $U_{ijkl}$,
\item[G3)] $i(A_{jk}+A_{kl} + A_{lj})=- d \log g_{jkl}$ on
 $U_{jkl}$,
\item[G4)] $(F_k-F_j)=dA_{jk}$ on $U_{jk}$.
\end{itemize}
G2) is
called the cocycle condition. The curvature of this gerbe-connection is the
globally-defined $3$-form $G\in \Lambda ^3(M)$, given locally by
\[
G=dF_i.
\]

\begin{figure}[ht!]
\centerline{\relabelbox\small 
\epsfysize 8cm
\epsfbox{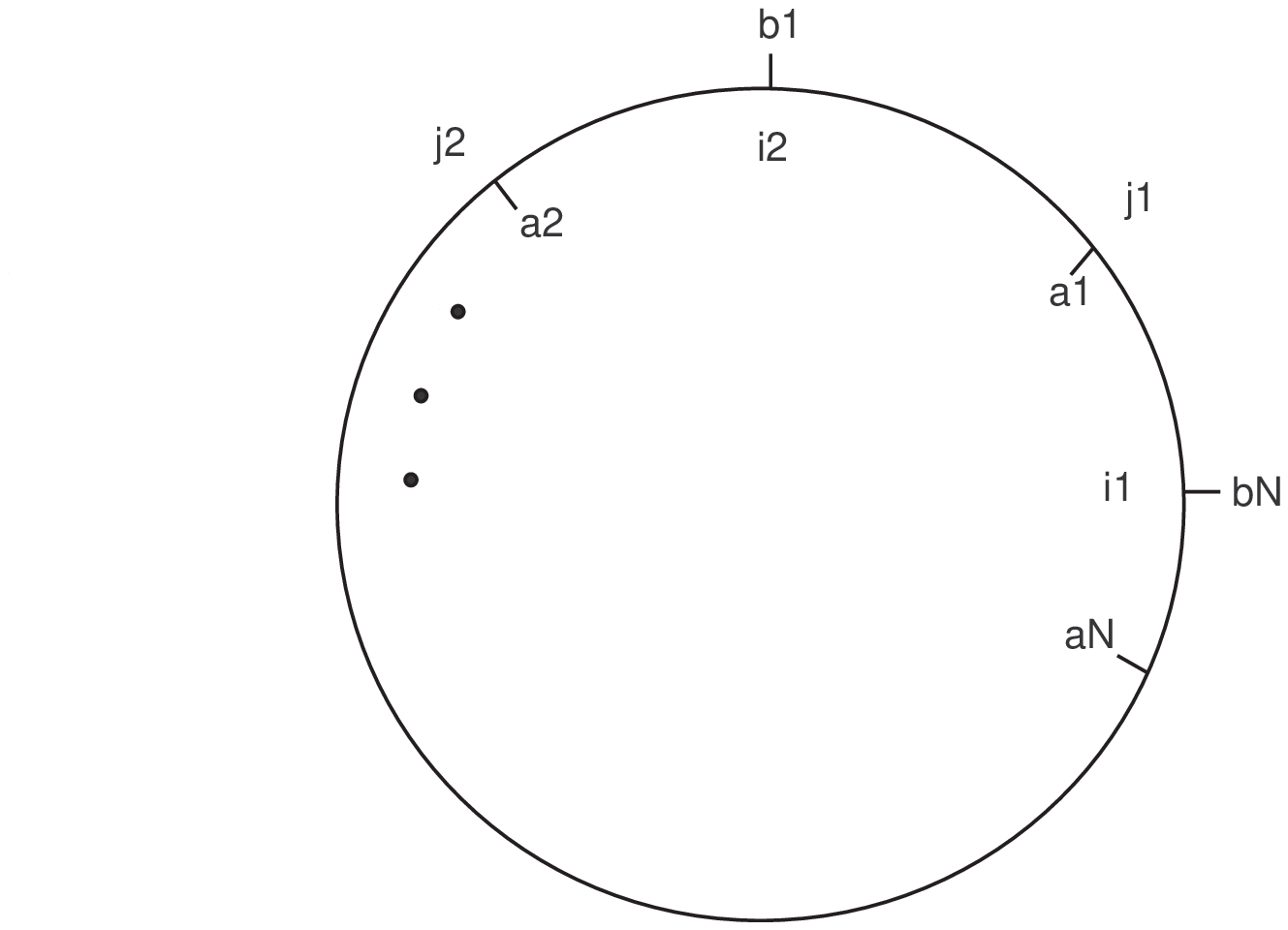} 
\relabel {a1}{$a_1$}
\relabel {a2}{$a_2$}
\relabel {aN}{$a_N$}
\relabel {b1}{$a'_1$}
\relabel {bN}{$a'_N$}
\relabel {i1}{$i_1$}
\relabel {i2}{$i_2$}
\relabel {j1}{$i'_1$}
\relabel {j2}{$ i'_2$}
\endrelabelbox }
\caption{Partitions for $\exp i \int_{(\ell,T,T')} (g, A)$}
\label{lTT}
\end{figure}

\begin{figure}[ht!]
\centerline{\relabelbox\small 
\epsfysize 3.5cm
\epsfbox{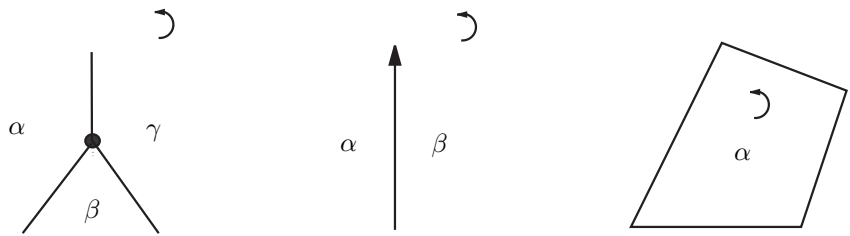} 
\relabel {a}{$\alpha$}
\relabel {b}{$\beta$}
\relabel {c}{$\gamma$}
\relabel {a1}{$\alpha$}
\relabel {b1}{$\beta$}
\relabel {a2}{$\alpha$}
\endrelabelbox }
\caption{Orientations for $v_{\alpha\beta\gamma}$, $e_{\alpha\beta}$ and
$f_{\alpha}$}
\label{orient}
\end{figure}

In the same way as before, we have
\begin{Thm}
A $U(1)$-gerbe with connection on $M$ gives rise to a  rank-$1$, embedded,
$2$-dimensional TQFT with target $M$ via
\begin{equation}
Z'(\ell ^\pm, T, T') = \exp \pm i \int_{(\ell,T,T')} (g, A)
\label{Z'l}
\end{equation}
\begin{equation}
Z(X,T)= \exp i \int_{(X,T)} (g,A,F),
\label{ZXT}
\end{equation}
where
\[
\begin{array}{lll}
\exp i \int_{(\ell,T,T')} (g, A) & = & g_{i_1i'_1i_2}(\ell({a}_1))
\exp i\int_{{a} _1}^{{{a}}'_1} \ell ^\ast(A_{i_2i'_1}).
\\ [15pt]
& & g_{i_2i'_1i'_2}(\ell({{a}}'_1))
\exp i\int_{{{a}}'_1}^{{{a}}_2} \ell ^\ast(A_{i_2i'_2}) \dots
\\ [15pt]
& & g_{i'_Ni'_1i_1}(\ell({{a}}'_N))
\exp i\int_{{{a}}'_N}^{{{a}}_1} \ell ^\ast(A_{i_1i'_1}),
\end{array}
\]
(see Figure \ref{lTT}, where we assume that the partitions can be put in
the general position indicated by subdividing and shifting regions), and
\[
\begin{array}{lll}
\exp i \int_{(X,T)} (g,A,F) & = & \prod^{int}_{
v_{\alpha\beta\gamma}} g_{i_\alpha i_\beta
i_\gamma}(X(v_{\alpha\beta\gamma})). \\ [15pt]
& & \prod^{int}_{
e_{\alpha\beta}}  \exp i \int_{e_{\alpha\beta}} X^\ast(A_{i_\alpha
i_\beta}).\\ [15pt]
& &
\prod^{int}_{
f_{\alpha}}  \exp i \int_{f_{\alpha}} X^\ast(F_{i_\alpha}),
\end{array}
\]
where the superscript on the products denotes that they are over internal
vertices, edges and faces, respectively  (see Figure \ref{orient} for the
orientations).
\label{thmgerbetqft}
\end{Thm}

\noindent\proof i) and ii) will
be shown using a Lemma to follow. iii) is immediate, since the factor
$Z'(X'|_C, T'_L, T'_R) $ gives precisely the missing internal contributions
in $Z(X',T')$ which come from the vertices and edges along $C$, as the
contributions from the transversal edges and the faces adjacent to $C$
vanish due to the ``constant at the boundary'' condition for $(X,T)$, and 
iv) is also immediate
 since there are no missing internal contributions in $Z(X',T')$, because the partitions on both sides of the curve are the same.

Indeed, to help prove the other two properties it is useful to observe
\[
\exp i \int_{(\ell , T, T')} (g,A) = \exp i \int_{(L , T\cup T')} (g,A, F)
\]
where $L:S^1 \times I\rightarrow M$ is given by $L({a}, x) = \ell
({a})$ and the partition of the annulus $T\cup T'$ is depicted in Figure
\ref{annulus}. This statement corresponds to c) of Proposition \ref{prop2} above.

We now introduce
 $3$-dimensional objects $(H,T)$, where $H$ is a smooth map
from an oriented $3$-manifold $V$ with or without boundary to $M$, and $T$ is a
labelled partition of $V$ into
regions $R_\alpha$, such that $H(R_\alpha)
\subset U_{i_\alpha}$,
by means of an embedded $4$-valent $2$-graph, which at each internal vertex
has four adjacent regions separated by six incident faces and four incident
edges. Each internal edge has three incident faces.  When $V$ has a
boundary, we will denote by $(\partial H, \partial T)$ the $2$-dimensional
object obtained by restricting $H$ and the labelled partition $T$ to the
oriented boundary $\partial V$ of $V$.

The Lemma from the previous section, one dimension up, now reads as follows.
\begin{Lem}
For a $3$-dimensional object $(H,T)$, with $H:V\rightarrow M$
\[
\exp i \int_{(\partial H, \partial T)} (g,A,F) = \exp i \int_V H^\ast (G).
\]
\end{Lem}
\proof
The contributions from the external transition functions may be
cancelled by introducing factors
\[
\exp i \int_{e_{\alpha\beta\gamma}} H^\ast
 (A_{i_\alpha i_\beta} + A_{i_\beta i_\gamma} + A_{i_\alpha i_\gamma} )
\]
along each internal edge, using G3), since the product of transition
functions at each internal vertex is $1$ because of the cocycle condition
 G2).
Using Stokes' theorem and G4), the contributions of all integrals along
edges may be replaced by factors
\[
\exp i \int_{f_{\alpha\beta}} H^\ast(F_{i_\beta} - F_{i_\alpha})
\]
for each internal face. Again using Stokes' theorem, these face integrals
may be replaced by
$$
\prod_\alpha \exp i \int_{R_\alpha} H^\ast(G) = \exp i \int_V H^\ast(G).
\eqno{\qed}$$

 We now  use the Lemma to prove the rest of Theorem \ref{thmgerbetqft}.

\begin{figure}[ht!]
\centerline{\relabelbox\small 
\epsfysize 6.5cm
\epsfbox{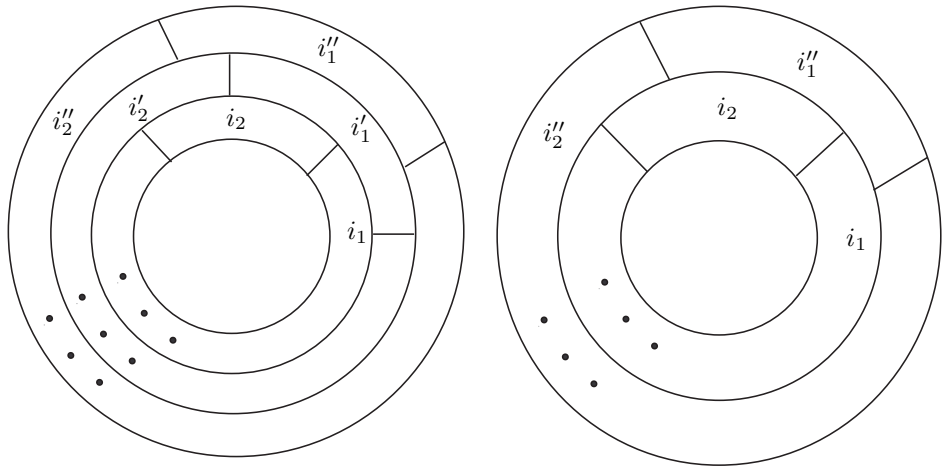} 
\relabel {i1}{$i_1$}
\relabel {i2}{$i_2$}
\relabel {j1}{$i'_1$}
\relabel {j2}{$i'_2$}
\relabel {k1}{$i''_1$}
\relabel {k2}{$i''_2$}
\relabel {l1}{$i_1$}
\relabel {l2}{$i_2$}
\relabel {m1}{$i''_1$}
\relabel {m2}{$i''_2$}
\endrelabelbox }
\caption{Partition on the top and bottom boundaries}
\label{annuli}
\end{figure}

 \proof[Proof of Theorem \ref{thmgerbetqft} (continued)] 
i) Using the remark after the proof of iii),
apply the Lemma to the map $H$ from the cylindrical shell $S^1\times
[0,1]^2$ to $M$, given by $H({a}, x,y) = \ell ({a})$, with the
partition on the top and bottom boundaries indicated in Figure
\ref{annuli}, and a vertical partition with constant labelling on the
inside and outside of the shell. It is easy to obtain an appropriate
partition of the cylindrical shell itself, by making a v-shaped ditch below
the annulus corresponding to $T'$ in the top boundary, and connecting the
circular bottom of the ditch via vertical surfaces to the middle circle of
the annulus in the bottom boundary.

ii) Apply the Lemma to a homotopy $H$ from $X$ to $X'$, whose image is
contained in the joint image of $X$ and $X'$. The boundary circle
contributions in ii) come from the annuli connecting the boundary circles
of $S$ and $S'$. It is always possible to obtain an appropriate partition
of the domain of $H$, since the partitions $T$ and $T'$ can be obtained
from each other by a sequence of elementary moves on labelled partitions of
surfaces, namely subdividing  a labelled face into faces with the same
label or recombining faces with the same label, and changing the label of a
face. It is simple to construct labelled volume partitions corresponding to
these moves. \endproof

Again the gerbe-induced embedded TQFT's have special properties, which are
proved in analogous fashion to Theorem \ref{thinthm}.
\begin{Thm} The rank-$1$, embedded, $2$-dimensional TQFT's of
Theorem {\em \ref{thmgerbetqft}} satisfy
\begin{itemize}
\item[\rm 1)] thin
invariance: when $(X,T)$ and $(X',T')$ have the same labelled
partition at their
boundaries and there is a relative homotopy $(H, \tilde{T})$ from $(X,T)$ to
$(X',T')$ satisfying ${\rm rank}\, DH \leq 2$, then
\[
Z(X,T)=Z(X',T'),
\]
\item[\rm 2)] smoothness: for any smooth $k$-dimensional family of objects
$(X(u),\!T(u)),$ $ u\in U\subset {\R}^k$, $Z(X(u),T(u))$ depends smoothly
on $u$.
\end{itemize}
\end{Thm}

We will conclude this section with the analogous theorem to
Theorem \ref{btqftcorr}.

\begin{figure}[ht!]
\centerline{\relabelbox\small 
\epsfysize 4cm
\epsfbox{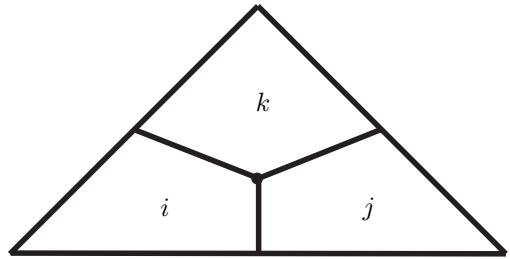} 
\relabel {i}{$i$}
\relabel {j}{$j$}
\relabel {k}{$k$}
\endrelabelbox }
\caption{Labelled partition of $\Delta$ for the gerbe transition function}
\label{gijk}
\end{figure}

\begin{figure}[ht!]
\centerline{\relabelbox\small 
\epsfysize 2.5cm
\epsfbox{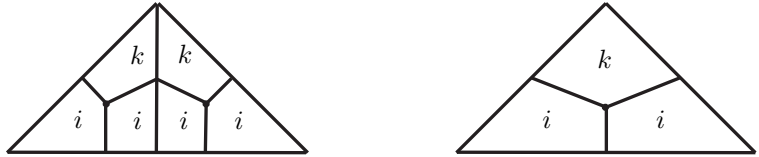} 
\relabel {i}{$i$}
\relabel {i1}{$i$}
\relabel {i2}{$i$}
\relabel {i3}{$i$}
\relabel {i4}{$i$}
\relabel {i5}{$i$}
\relabel {k}{$k$}
\relabel {k1}{$k$}
\relabel {k2}{$k$}
\endrelabelbox }
\caption{Proof of G1)}
\label{G1pf}
\end{figure}

\begin{figure}[ht!]
\centerline{\relabelbox\small 
\epsfysize 5cm
\epsfbox{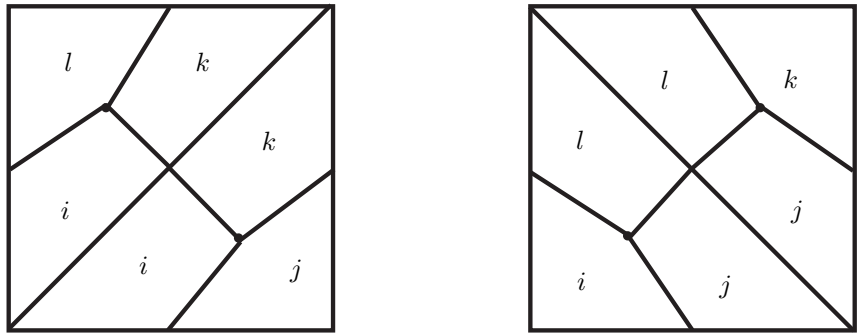} 
\relabel {i}{$i$}
\relabel {i1}{$i$}
\relabel {i2}{$i$}
\relabel {j}{$j$}
\relabel {j1}{$j$}
\relabel {j2}{$j$}
\relabel {k}{$k$}
\relabel {k1}{$k$}
\relabel {k2}{$k$}
\relabel {l}{$l$}
\relabel {l1}{$l$}
\relabel {l2}{$l$}
\endrelabelbox }
\caption{Proof of G2)}
\label{G2pf}
\end{figure}

\begin{Thm}
There is a one-to-one correspondence between gerbes with connection on $M$
and smooth, thin invariant, rank-$1$, embedded, $2$-dimensional TQFT's with
target $M$.
\label{gtqftcorr}
\end{Thm}
\proof
We have already seen the correspondence in one direction.
Suppose we are given an embedded TQFT with the stated properties. We may
reconstruct a gerbe with connection on $M$ as follows.
\begin{equation}
g_{ijk}(y) = Z(\Delta_y, ijk),
\label{gfromZ2}
\end{equation}
where $\Delta_y$ is the constant map from a standard $2$-simplex $\Delta$ to $y$, and
$ijk$ denotes the labelled partition obtained from connecting the midpoint of each edge of $\Delta$ to a trivalent vertex at the centre of $\Delta$, and assigning $i,j$ and $k$ to the three regions of this partition (Figure \ref{gijk}).
Properties G1) and G2) now follow from applying $Z$ to the partial gluing morphisms and isomorphisms depicted in Figures \ref{G1pf} and \ref {G2pf}. Note that we need partial gluing morphisms here, and not gluing morphisms, since only a part of the boundary of the simplices is glued.

Next, the connection $1$-forms are defined as follows:
\begin{equation}
i(A_{jk})_y(v)= \frac{d}{dt} \log Z(Q_{t}, jk) |_{t=0}
\label{AfromZ2}
\end{equation}
where we take a smooth path $q:(-\epsilon , \epsilon)\rightarrow U_{jk}$ such that
$q(0)=y, \, \dot{q}(0)=v$, and $Q_{t}$ is a map from $[0,1]_u\times [0,1]_s$ to $U_{jk}$ which does not depend on $u$, and as a function of $s$ is a path from $y$ to $q(t)$ along $q$, reparametrized to be constant at the endpoints.
Furthermore, $jk$ denotes the labelled partition which assigns $j$ to the region $u\in [0,1/2]$ and $k$ to the region $u\in [1/2,1]$. $A_{jk}$ is well-defined by an identical argument to that used to show that $A_j$ is well-defined in the proof of Theorem \ref{btqftcorr}. In particular we have
\[
\frac{d}{dt} \log Z(Q_{t}, jk) |_{t=0} = - \frac{d}{dt} \log Z(R_{Q(t)}, jk) |_{t=0}= -(dh)_y(v),
\]
where $R_x: [0,1]_u\times [0,1]_s\rightarrow M$ is defined in terms of the radial paths of that proof, namely $R_x(u,s)=r_x(s)$, and the function $h$ is given by $h(x)=\log Z(R_x, jk)$ in a neighbourhood of $y$.

\begin{figure}[ht!]
\centerline{\relabelbox\small 
\epsfysize 3.5cm
\epsfbox{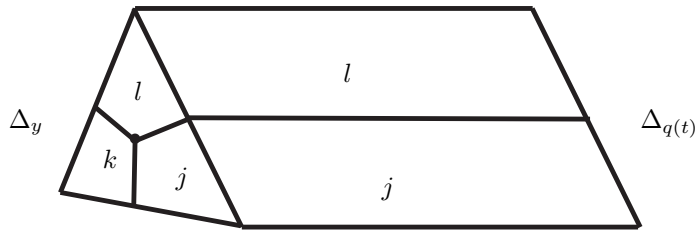} 
\relabel {i}{$j$}
\relabel {i1}{$j$}
\relabel {j}{$k$}
\relabel {k}{$l$}
\relabel {k1}{$l$}
\relabel {Dy}{$\Delta_y$}
\relabel {Dq}{$\Delta_{q(t)}$}
\endrelabelbox }
\caption{Labelled partition for the proof of G3)}
\label{prism}
\end{figure}

To show that the transition functions and connection $1$-forms satisfy G3), consider the function
$\varphi:  (-\epsilon, \epsilon) \rightarrow i {\R}$ defined by
\[
\varphi(t) = \log Z(P_t, jkl),
\]
where $P_t$ is the map from the prism surface to $U_{ijk}$ given by $\Delta_y$ and $\Delta_{q(t)}$ on the two triangular faces, and $Q_{t}$ on the three square faces, and the labelled partition $ijk$ is indicated in Figure \ref{prism}. By the thin invariance of $Z$, $\varphi$ is constant, since there are natural thin homotopies between $P_t$'s for different values of $t$ (the image of these homotopies is contained in the image of $q$). Thus
\[
\begin{array}{lll}
0 & =  & \frac{d}{dt} [ -\log Z(\Delta_y, jkl) + \log Z(\Delta_{q(t)}, jkl)\, +
\\
& & \log (Z(Q_{t}, jk) + \log (Z(Q_{t}, kl) +
\log (Z(Q_{t}, lj) ] |_{t=0}
\\
& = & i(A_{jk} + A_{kl} + A_{lj})_y(v) + (d \log g_{jkl})_y(v).
\end{array}
\]

Finally, the connection $2$-forms are defined as follows:
\begin{equation}
i(F_j)_y(v,w) = \frac{\partial^2}{\partial t \partial u}
\log Z(Q_{t,u},j) |_{(t,u)=(0,0)},
\label{FfromZ}
\end{equation}
where we take a smooth map $Q: (-\epsilon , \epsilon)_r \times
(-\epsilon , \epsilon)_s
\rightarrow U_{i}$ such that
\[
Q(0,0)=y,\qquad \frac{\partial Q}{\partial r}(0,0)=v,\qquad \frac{\partial Q}{\partial s}(0,0)=w,
\]
and set $ Q_{t,u}$ to be a $2$-path (map from the Cartesian product of two intervals to $M$) from $Q(0,.)$ to $Q(t,.)$  and from $Q(.,0)$ to $Q(.,u)$  along $Q$, reparametrized to be constant on the boundary. Furthermore, $j$ denotes the labelled partition which assigns $j$ to the whole domain of 
$Q_{t,u}$. 

We sketch the proof that $F_j$ is well-defined. First, given $Q$, $F_j$ does not depend on the choice of $Q_{t,u}$ because of the thin invariance of $Z$. Thus it remains to show that the definition does not depend on the choice of $Q$. We construct a second family of maps $\tilde{Q}_{t,u}:D\rightarrow M$, where $D$ denotes the unit disk, such that $\tilde{Q}_{t,u}(1,0)=y$, $ \tilde{Q}_{t,u}$ restricted to the four quadrants of $\partial D$ coincides with ${Q}_{t,u}$ restricted to the four edges of the square, in an obvious sense, and $\tilde{Q}_{t,u}$ restricted to any chord in $D$ ending at $(1,0)$ is the radial path in $M$ (as in the proof of Theorem \ref{btqftcorr})  from the image of the starting point of the chord to $y$. Combining the maps ${Q}_{t,u}$ and $\tilde{Q}_{t,u}$ gives a $2$-parameter family of maps from $S^2$ to $M$, equal to the trivial map at $t=u=0$. Applying the higher version of Barrett's lemma proved in \cite{MP}, we have
\[
\frac{\partial^2}{\partial t \partial u}\log Z(Q_{t,u},j) |_{(t,u)=(0,0)} = 
\frac{\partial^2}{\partial t \partial u}\log Z(\tilde{Q}_{t,u},j) |_{(t,u)=(0,0)}. 
\]
The right-hand-side can be shown to be equal to $(dB)_y(v,w)$, where $B$ is a smooth $1$-form defined in a neighbourhood of $y$, and is thus independent of the choice of $Q$.

\begin{figure}[ht!]
\centerline{\relabelbox\small 
\epsfysize 5cm
\epsfbox{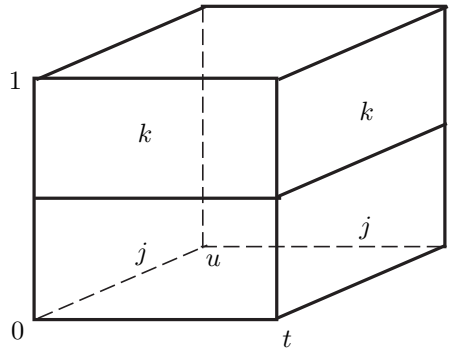} 
\relabel {0}{$0$}
\relabel {1}{$1$}
\relabel {t}{$t$}
\relabel {u}{$u$}
\relabel {i}{$j$}
\relabel {i1}{$j$}
\relabel {j}{$k$}
\relabel {j1}{$k$}
\endrelabelbox }
\caption{Labelled partition for the proof of G4)}
\label{G4pf}
\end{figure}

To show that equation G4) holds for the connection $1$- and $2$-forms, consider the function 
\[
\Phi: \, (-\epsilon, \epsilon) \,\times\, (-\epsilon,\epsilon)
\rightarrow i {\R}, \qquad
\Phi(t,u) = \log Z(B_{t,u}, jk),
\]
where $B_{t,u}$ maps from the boundary of the cube 
$[0,t]_r\times [0,u]_s \times [0,1]_p$ to $U_{ij}$, and is given by 
$ Q_{t,u}$ on $p=0,1$, and is constant in $p$ on the sides $r=0,1,\, s=0,1$. Furthermore, $jk$ denotes the labelled partition indicated in Figure \ref{G4pf}. $\Phi$ is constant, by the thin invariance of $Z$, since the images for all arguments are contained in the image of $Q$, and G4) follows from writing out the equation $\partial^2 \Phi/\partial t \partial u (0,0) =0$. 

Finally we show the correspondence is one-to-one. Let 
$(g_{ijk}, A_{jk}, F_k)_{i,j,k\in J}$ be a gerbe with connection, and let $(Z',Z)$ be the TQFT which arises from it (Theorem \ref{thmgerbetqft}). Now take 
$(\tilde{g}_{ijk},\tilde{A}_{jk},\tilde{F}_k )_{i,j,k\in J}$  to
be the gerbe and connection reconstructed from $(Z',Z)$ via equations (\ref{gfromZ2}), (\ref{AfromZ2}) and (\ref{FfromZ}). We have
\[
\tilde{g}_{ijk}(y)= Z(\Delta_y, ijk) =g_{ijk}(y)
\]
using equation (\ref{ZXT}) in the second equality, since the edge and face integrals on the right-hand-side don't contribute, as the map $\Delta_y$ is constant. Next we have
\begin{eqnarray*}
i(\tilde{A}_{jk})_y(v) &=&\frac{d}{dt} \log Z(Q_{t}, jk) |_{t=0}=
i \frac{d}{dt} \int_0^1 Q_t^*(A_{jk}) |_{t=0} \\
&=& i \frac{d}{dt}  \int_0^t q^*(A_{jk}) |_{t=0} = i(A_{jk})_y(v)
\end{eqnarray*}
where, in the first equality, we choose 
$Q_t: [0,1]^2\rightarrow M$ defined by $Q_t(u,s)=q(st)$, reparametrized to be constant at the boundary, in the second equality we use equation (\ref{ZXT}), and the fact that there are no vertices and the face integrals don't contribute as $Q_t$ is constant in $u$, and, in the third equality, we use the reparametrization invariance of the integral and make a substitution. Thirdly, we have
\begin{eqnarray*}
i(\tilde{F}_k)_y(v,w)& = & \frac{\partial^2}{\partial t \partial u}
\log Z(Q_{t,u},k) |_{(t,u)=(0,0)} \\
& = & i \frac{\partial^2}{\partial t \partial u}
 \int_{[0,1]^2} Q_{t,u}^*(F_k) |_{(t,u)=(0,0)} \\
& = & i \frac{\partial^2}{\partial t \partial u}
 \int_0^t \int_0^u Q^*(F_k) |_{(t,u)=(0,0)}=  ({F}_k)_y(v,w)
\end{eqnarray*}
where in the first equality we choose $Q_{t,u}: [0,1]^2\rightarrow M$ defined by $Q_{t,u}(r,s)=Q(rt,su)$, reparametrized to be constant at the boundary, in the second equality we use equation  (\ref{ZXT}), and in the third equality we use the reparametrization invariance of the integral and make two substitutions.

Conversely, let $(Z',Z)$ be a TQFT, as specified, and let 
$(g_{ijk}, A_{jk}, F_k)_{i,j,k\in J}$ be the gerbe with connection reconstructed from it via equations (\ref{gfromZ2}), (\ref{AfromZ2}) and (\ref{FfromZ}). 
Now take $(\tilde{Z}', \tilde{Z})$ to be the TQFT that this bundle with connection gives rise to (Theorem \ref{thmgerbetqft}). We have
\begin{eqnarray*}
\tilde{Z}'(\ell^+,T,T') & = & g_{i_1i'_1i_2}(\ell({a}_1))
\exp i\int_{{a} _1}^{{{a}}'_1} \ell ^\ast(A_{i_2i'_1}) \dots \\
&=& Z(\Delta_{\ell({a}_1)}, i_1i'_1i_2) Z(L_{11'}, i_2i'_1) \dots \\
&=& Z(L,T\cup T') = Z'(\ell^+,T,T').
\end{eqnarray*}
Here in the second equality we use equation (\ref{gfromZ2}) and the relation
\[
\exp i\int_{{a} _1}^{{{a}}'_1} \ell ^\ast(A_{i_2i'_1}) 
= Z(L_{11'}, i_2i'_1)
\]
where $L:S^1\times I \rightarrow M$ is given by $L({a}, x) = \ell({a})$, $L_{11'}$ denotes $L$ restricted to the interval $[{a}_1,{a}'_1]$, and $i_2i'_1$ denotes the partition with $i_2$ on the inside and $i'_1$ on the outside of the strip (see Figures \ref{lTT} and \ref{annulus}). This relation is shown in analogous fashion to equation (\ref{eintA=Z}) in the proof of Theorem 
\ref{btqftcorr}. In the third equality we use the partial gluing property iv) 
of Definition \ref{def2dTQFT}
to put all the pieces together, as well as an isomorphism between $L$ with the partition it acquires in this way, and $L$ with the partition $T\cup T'$ of Figure \ref{annulus}. The last equality is Proposition \ref{prop2} c). 
Next we have
\[
\tilde{Z}'(\ell^-,T,T')= \tilde{Z}'(\ell^+,T,T')^{-1}
= {Z}'(\ell^+,T,T')^{-1} = {Z}'(\ell^-,T,T'),
\]
using Proposition \ref{prop2} d) in the last equality. Finally we have
\[
\begin{array}{lll}
\tilde{Z}(X,T) & = & \exp i \int_{(X,T)} (g,A,F) \\ [15pt]
& = & \prod^{int}_{v_{\alpha\beta\gamma}}
Z(\Delta_{X(v_{\alpha\beta\gamma})},i_\alpha i_\beta i_\gamma). \\ [15pt]
& & \prod^{int}_{e_{\alpha\beta}}
  Z(\epsilon_{X,\alpha\beta},i_\alpha i_\beta).\\ [15pt]
& &
\prod^{int}_{f_{\alpha}}  
Z(\phi_{X,\alpha},i_\alpha) = Z(X,T)
\end{array}
\]
where we have used equation (\ref{gfromZ2}) for the vertex factors, we have introduced the $2$-path $\epsilon_{X,\alpha\beta}:[0,1]_s\times [0,1]_t \rightarrow M$, constant in $t$ and equal to the path $X$ restricted to the edge $e_{\alpha\beta}$ as a function of $s$, and the labelling $i_\alpha i_\beta$ which denotes assigning $i_\alpha$ to the region $t\geq 1/2$ and $i_\beta$ to the region $t\leq 1/2$. We have also introduced $\phi_{X,\alpha}:f_\alpha \rightarrow M$ being $X$ restricted to $f_\alpha$, and the labelling $i_\alpha$ which denotes assigning $i_\alpha$ to the whole of $f_\alpha$. The equality for the edge factors is shown in analogous fashion to equation (\ref{eintA=Z}) in the proof of Theorem \ref{btqftcorr} and the equality for the face factors is shown by an analogous argument for $2$-forms. The final equality uses iv) 
of Definition \ref{def2dTQFT} to put all the pieces together, as well as an isomorphism between the $2$-dimensional object thus obtained and $(X,T)$. \endproof

\section{Comments}
\label{Comm}

In this direct approach using the local functions and forms which make up
bundles or gerbes with connection, we were able to avoid the use of
trivializations of the pull-back bundle or gerbe, which appeared in various
proofs in \cite{MP}. Indeed these trivializations may not exist if the
topology of the object's domain is non-trivial. We note that the integral
formulae given here do not require $\pi_1(M)$ to be trivial, whereas
non-simply-connected manifolds $M$ created some technical hurdles in
\cite{MP}.

Bundles, gerbes and higher generalizations sit nicely inside a variety of
dimensional ladders, and several people have
conjectured links between this chain of higher-order geometries and
state-sum models of TQFT. The constructions and methods of proof here
have a distinct state-sum flavour, and there
are intriguing suggestions of a graded
integration theory, involving objects of all dimensions up to the dimension
being integrated. This could be useful for understanding
broader dimensional ladders for TQFT \cite{BaeDol}, which normally involves
just the top dimension and one dimension lower, or TQFT with corners. It
would also be interesting  to see whether there are links between these
constructions and the state-sum models for quantum gravity proposed by
Barrett and Crane \cite{BarCra} and Mikovic \cite{Mik}.

Approaches to non-abelian gerbes with connection have been
studied by various authors \cite{BreMes, Att, Bae}. We could have obtained
solutions of the higher-rank $1$-dimensional embedded TQFT equations,
mentioned after Definition \ref{def1dTQFT}, from
nonabelian bundles with connection, via path-ordered exponentials. The
geom\-et\-ric-combinatorial approach in this article may suggest non-abelian
solutions for the higher-rank $2$-dimensional embedded TQFT equations,
which would then be strong candidates to be called non-abelian gerbes with
connection.

\Addressesr


\begin{thebibliography}



\bibitem{Ati89}
{\bf M Atiyah},
{\it Topological Quantum Field Theories},
 Publ. Math. Inst. Hautes Etudes Sci. 68 (1989) 175--186


\bibitem{Att}
{\bf R Attal},
{\it Combinatorics of Non-Abelian Gerbes with Connection and Curvature},
\arxiv{math-ph/0203056}



\bibitem{Bae}
{\bf J Baez},
  {\it Higher Yang-Mills Theory},
  \arxiv{hep-th/0206130}


\bibitem{BaeDol}
{\bf J Baez, J Dolan},
 {\it Higher-dimensional algebra and topological quantum field theory},
  J. Math. Phys. 36(11) (1995) 6073--6105

\bibitem{Bar91}
{\bf J\,W Barrett},
{\it  Holonomy and path structures in general relativity and {Y}ang-{M}ills
  theory},
Int. J. Theor. Phys. 30(9) (1991) 1171--1215


\bibitem{BarCra}
{\bf J\,W Barrett, L Crane},  
{\it Relativistic spin networks and quantum gravity},
J.Math.Phys. 39 (1998) 3296--3302


\bibitem{BreMes}
{\bf L Breen, W Messing},
 { \it Differential Geometry of Gerbes}, preprint (2001)
  \arxiv{math.AG/0106083}


\bibitem{BriTur00}
{\bf M Brightwell, P Turner},
  {\it Representations of the homotopy surface category of a simply connected space},
  J. Knot Theory Ramifications 9 (2000) 855--864 

\bibitem{Bry}
{\bf J-W Brylinski},
  {\it Loop spaces, characteristic classes and geometric quantization},
  volume 107 of {\it Progress in Mathematics},
  Birkhauser (1993)


\bibitem{BunTurWil}
{\bf U Bunke, P Turner, S Willerton},
  {\it Gerbes and homotopy quantum field theories},
  \arxiv{math.AT/0201116}


\bibitem{CP94}
{\bf A Caetano, R Picken},
 {\it An axiomatic definition of holonomy},
 Int. J. Math. 5(6) (1994) 835--848 

\bibitem{gaw}
{\bf K.~Gaw\c{e}dzki},
{\it  Topological actions in two-dimensional quantum field theory}, 
  In {\it Non-perturbative quantum field theory}, G. 't Hooft et al eds., Plenum Press (1988)

\bibitem{gaw:rei} 
{\bf K Gaw\c{e}dzki, N Reis},
{\it  WZW branes and gerbes},
  Rev. Math. Phys. 14 (2002) 1281--1334


\bibitem{Gir}
{\bf J Giraud},
  {\it Cohomologie non-abelienne}, volume 179 of {\it Grundl.},
  {S}pringer-{V}erlag (1971)


\bibitem{Hi99}
{\bf N Hitchin},
{\it  Lectures on special {L}angrangian submanifolds},
  in {\it {W}inter {S}chool on {M}irror {S}ymmetry, {V}ector {B}undles
and {L}agrangian {S}ubmanifolds (Cambridge, MA, 1999)}, AMS/IP Stud. Adv. Math. 23, American Mathematical Society (2001) 151--182

\bibitem{MP}
{\bf M Mackaay, R Picken},
  {\it Holonomy and Parallel Transport for Abelian Gerbes},
Adv. Math. 170 (2002) 287--339 

\bibitem{Mik}
{\bf A Mikovic},
  {\it Spin Foam Models of Matter Coupled to Gravity},
Class. Quant. Grav. 19 (2002) 2335--2354 


\bibitem{PS}
{\bf R Picken, P Semi\~{a}o},
  {\it A classical approach to TQFT's}, preprint (2002)
  \arxiv{math.QA/0212310}


\bibitem{Rod01}
{\bf G Rodrigues},
  {\it Homotopy quantum field theories and the homotopy cobordism category
  in dimension $1+1$}, 
J. Knot Theory Ramifications 12(3) (2003) 287--319 


\bibitem{Seg68}
{\bf G Segal},
  {\it Classifying spaces and spectral sequences},
  Inst. Hautes \'{E}tudes Sci. Publ. Math. 34 (1968) 105--112


\bibitem{Seg89}
{\bf G Segal},
  {\it Two-Dimensional Conformal Field Theories and Modular Functors},
  in {\it Proceedings of the XIth International Congress on
  Mathematical Physics, Swansea, 1988}, Adam Hilger (1989) 22--37 

\bibitem{Seg01}
{\bf G Segal},
  {\it Topological Structures in String Theory},
Phil. Trans. R. Soc. Lond. A 359 (2001) 1389--1398 

\bibitem{SegCFT}
{\bf G Segal},
  {\it The definition of  Conformal Field Theory},
  Cambridge Univ. Press, to appear.


\bibitem{Tur99}
{\bf V Turaev},
  {\it Homotopy field theory in dimension 2 and group-algebras}, e-print
  \arxiv{math.QA/9910010}

\bibitem{Turn04}
{\bf P Turner},
  {\it A functorial approach to differential characters},
Alg. Geom. Topol. 4(6) (2004) 81--93 


\bibitem{Wit88}
{\bf E Witten},
  {\it Topological Quantum Field Theory},
Commun. Math. Phys. 117 (1989) 353--386 


\end{thebibliography}
\end{document}